\documentclass[10pt]{article}

\usepackage{bm}
\usepackage{amsmath}
\usepackage{amsfonts}
\usepackage{amscd}
\usepackage{hyperref}
\usepackage{pstricks}
\usepackage{latexsym}
\usepackage{enumerate}
\usepackage{epsfig}
\usepackage{subfigure}
\usepackage{graphicx}
\usepackage{caption}
\usepackage{subcaption}
\usepackage{framed,fancybox}
\usepackage{bbm}
\usepackage{multirow}
\usepackage[T1]{fontenc}
\usepackage{threeparttable}
\usepackage{rotating}
\usepackage{color,soul}
\usepackage{physics}
\usepackage[version=4]{mhchem}
\usepackage{siunitx}
\usepackage{cite}
\usepackage{footnpag}			      	% make footnote symbols restart on each page
\usepackage{longtable,tabularx}
\newcommand{\m}[1]{{\bf{#1}}}

\newcommand{\C}[1]{{\cal {#1}}} %script letters

\def\BibTeX{{\rm B\kern-.05em{\sc i\kern-.025em b}\kern-.08em
    T\kern-.1667em\lower.7ex\hbox{E}\kern-.125emX}}
\begin{document}

\title{\bf Desensitized Optimal Guidance Using \\ Adaptive Radau Collocation\thanks{The authors gratefully acknowledge support for this research from the U.S. National Science Foundation under grant CMMI-2031213, from the U.S. Office of Naval Research under grant N00014-22-1-2397, and from the U.S. Air Force Research Laboratory under grant FA8651-21-F-1041.}
}
\author{Katrina L.~Winkler${}^{*}$\thanks{${}^{*}$Ph.D.~Student, Department of Mechanical and Aerospace Engineering. E-mail: k.winkler@ufl.edu.} and Anil V.~Rao${}^{\dag}$\thanks{${}^{\dag}$Professor, Department of Mechanical and Aerospace Engineering. E-mail: anilvrao@ufl.edu.} \\ \\ {\em University of Florida} \\ {\em Gainesville, FL 32611}}

\maketitle

\begin{abstract}
An optimal guidance method is developed that reduces sensitivity to parameters in the dynamic model. The method combines a previously developed method for guidance and control using adaptive Legendre-Gauss-Radau (LGR) collocation and a previously developed approach for desensitized optimal control.  Guidance updates are performed such that the desensitized optimal control problem is re-solved on the remaining horizon at the start of each guidance cycle.  The effectiveness of the method is demonstrated on a simple example using Monte Carlo simulation.  It is found that the method reduces variations in the terminal state as compared to either desensitized optimal control without guidance updates or a previously developed method for optimal guidance and control.
\end{abstract}

\section{Introduction}
The goal of an optimal control problem is to determine the state and control of a controlled dynamical system that optimizes a specified performance index while satisfying dynamic constraints, path constraints, and boundary conditions. Due to their complexity, most optimal control problems cannot be solved analytically and, thus, must be solved numerically.  Numerical methods for solving optimal control problems fall into two categories: indirect methods and direct methods.  In an indirect method, calculus of variations is employed to formulate the first-order optimality conditions, leading to a Hamiltonian boundary-value problem (HBVP).  A numerical solution of the HBVP is then obtained. In a direct method, the state and control are parameterized and the optimal control problem is transcribed into a finite-dimensional nonlinear programming problem (NLP) which is solved using well known optimization methods \cite{betts2010practical,biegler2009large}.

Optimal control problems are typically formulated with a reference (nominal) dynamic model and the optimized trajectory and control are obtained in an open-loop manner along the horizon of interest using this reference model. In the absence of uncertainty, the reference optimal control obtained using the reference dynamic model would suffice.  In a real system, however, uncertainties will exist.  As a result, using the reference optimal control will lead to perturbations from the reference optimal trajectory.  In the case of sufficiently large uncertainties, employing the reference optimal control for any significant duration will lead to large perturbations in the state and, in order to reduce such perturbations, some form of correction will be required.  Often, course corrections in the form of guidance updates are performed where the optimal control problem may be solved periodically (that is, at specified guidance update times), thus leading to closed-loop optimal control. 

% The first is inner-loop guidance which implements feedback through introduction of control-like gain functions \cite{seywald2019desensitized}. The second method is known as outer-loop guidance which resolves the optimal control problem on the unexpired horizon to redetermine the optimal control from the current state to the final time \cite{Dennis_2019}. 

% Trajectory optimization uses optimal control to determine a reference state and control which are then tracked by guidance algorithms to ensure the vehicle at hand reaches its desired final state while satisfying constraints and minimizing the user-defined cost. Trajectory optimization of the nominal trajectory is an open loop process in which corrections in the optimal control are not performed. However, in real life applications, the system is subject to perturbations including but not limited to uncertainties in parameters that the system's dynamics may be a function of. Thus, guidance methods are required to allow for corrections in the control so that the vehicle's trajectory remains optimal despite state perturbations. Allowing for corrections in the control is known as closed-loop guidance and can be done one of two ways. The first is inner-loop guidance which implements feedback through introduction of control-like gain functions \cite{seywald2019desensitized}. The second method is known as outer-loop guidance which resolves the optimal control problem on the unexpired horizon to redetermine the optimal control from the current state to the final time \cite{Dennis_2019}. 

In the application of optimal control for high performance vehicles, the vehicle is subject to uncertainty in model parameters.  When designing a reference trajectory, it is desirable to reduce sensitivities to these parametric uncertainties to promote robustness while minimizing the error in the final state in response to perturbations in the state anywhere along the trajectory. The process of determining the state and control that reduces sensitivity to parametric uncertainty is known as desensitized optimal control. The study of desensitized optimal control first appeared in Ref.~\cite{seywald2003desensitized} for a simple optimal control problem. An extension of the work of Ref.~\cite{seywald2003desensitized} for problems with constraints can be found in  Ref.~\cite{seywald2003desensitizedConstraints}. Next, Ref.~\cite{shen2008desensitizing} and \cite{shen2010desensitizing} use desensitized optimal control to solve the Mars pinpoint landing problem. The work of Ref.~\cite{shen2008desensitizing} and \cite{shen2010desensitizing} was then extended in Ref.~\cite{xu2015robust} to study perturbations resulting from parametric uncertainties in aerodynamic characteristics and atmospheric density and was further implemented using direct collocation and nonlinear programming in Ref.~\cite{li2011mars}. Furthermore, Ref.~\cite{zimmer2005reducing} studied a relationship between covariance trajectory shaping and desensitized optimal control, while Ref.~\cite{small2010optimal} studied this same relationship in conjunction with trajectory design.  Next, desensitized trajectory optimization was studied in Ref.~\cite{Makkapati_2018} where the sensitivity dynamics were explored as functions of the partial derivatives of the original dynamics with respect to the state and parameters with uncertainties. In particular, a simplified form of the sensitivity dynamics from Ref.~\cite{seywald2003desensitized} was derived. This method was then applied to hypersonic trajectory optimization of a reentry problem \cite{Makkapati_2021} in which parameter uncertainties existed in the parasitic drag and scaling height. The cost was augmented to include the expected deviations in a user-defined penalty term.

The objective of this work is to develop a method for optimal guidance that desensitizes the reference control to parametric uncertainties while providing guidance updates to allow for corrections in the desensitized optimal control. The method developed employs a closed-loop adaptation of the method developed in Ref.~\cite{Makkapati_2021} while employing a guidance strategy that is based on the work of Ref.~\cite{Dennis_2019}. The adaptation of the method of Ref.~\cite{Makkapati_2021} offers an approach for  desensitized optimal control while the inclusion of the method developed in Ref.~\cite{Dennis_2019} provides an efficient approach for performing guidance updates using adaptive Legendre-Gauss-Radau (LGR) collocation.  The method developed in this paper is demonstrated on a simple example where it is shown that combining guidance updates via solving a desensitized optimal control problem reduces the error in the final state compared with using either open-loop desensitized optimal control or optimal guidance alone. % The optimal control software GPOPS-II was used to determine the DOC solution \cite{patterson2014gpops}. 

\section{Bolza Optimal Control Problem}

Without loss of generality, consider the following Bolza optimal control problem in terms of the elapsed time, $t$, that has the domain $t \in [t_0,t_f]$. Determine the state $\m{x}(t)$ and the control $\m{u}(t)$ that minimizes the cost functional
\begin{equation}
    J = \mathcal{M}(\m{x}(t_0),t_0,\m{x}(t_f),t_f) + \int_{t_0}^{t_f} \mathcal{L}(\m{x}(t),\m{u}(t),t)dt,
\end{equation}
subject to the dynamic constraints
\begin{equation}
    \dot{\m{x}} = \m{f}(\m{x}(t),\m{u}(t),t),
\end{equation}
boundary conditions
\begin{equation}
    \m{b}_{\min} \leq \m{b}(\m{x}(t_0),t_0,\m{x}(t_f),t_f) \leq \m{b}_{\max},
\end{equation}
and inequality path constraints
\begin{equation}
    \m{c}_{\min} \leq \m{c}(\m{x}(t),\m{u}(t),t) \leq \m{c}_{\max}.
\end{equation}

Suppose now that the trajectory time domain, $t\in[t_0,t_f]$, is mapped to a mesh time domain, $\tau \in [-1, +1]$ where
\begin{equation}\label{tau-to-t}
    t \equiv t(\tau,t_0,t_f) = \frac{t_f+t_0}{2}\tau + \frac{t_f-t_0}{2}. 
\end{equation}
Next, let the mesh domain be divided into $K$ mesh intervals $\{\C{I}_1,\ldots,\C{I}_K\}$ such that $\C{I}=[T_{k-1},T_k],\; (k=1,\ldots,K)$ where the mesh points $(T_0,\ldots,T_K)$ are defined such that $-1=T_0<T_1<\cdots<T_{K-1}<T_K=+1$. Furthermore, let $\m{x}^k(\tau)$ and $\m{u}^k(\tau)$ denote, respectively, the state and control in mesh interval $\C{I}_k,\; (k=1,\ldots,K)$.  The Bolza optimal control problem can then be written in terms of the mesh time domain as follows: Minimize the objective functional
\begin{equation}
    J = M(\m{x}^{(1)}(T_0),t_0,\m{x}^{(K)}(T_K),t_f) \\+ \frac{t_f - t_0}{2} \sum_{k=1}^K \int_{T_{k-1}}^{T_k} \mathcal{L} \left(\m{x}^{(k)}(\tau),\m{u}^{(k)}(\tau),t\right)d\tau,
\end{equation}
subject to the dynamic constraints
\begin{equation}
    \dot{\m{x}}^{(k)}(\tau) = \frac{t_f - t_0}{2} \m{f}(\m{x}^{(k)}(\tau),\m{u}^{(k)}(\tau),t),\; (k=1,\ldots,K), 
\end{equation}
the path constraints
\begin{equation}
    \m{c}_{\min} \leq \m{c}(\m{x}^{(k)}(\tau),\m{u}^{(k)}(\tau),t) \leq \m{c}_{\max},\; (k=1,\ldots,K),
\end{equation}
and the boundary conditions
\begin{equation}
    \m{b}_{\min} \leq \m{b}(\m{x}(T_0),t_0,\m{x}(T_K),t_f) \leq \m{b}_{\max}.
\end{equation}
Finally, continuity in the state at each interior mesh point $T_k,\; (k=1,\ldots,K-1)$ is enforced via the constraint $\m{x}^{k}(T_{k})=\m{x}^{k+1}(T_{k}),\; (k=1,\ldots,K-1)$.

\section{Legendre-Gauss-Radau Collocation Method}
 The proposed desensitized optimal control guidance scheme employs the $hp$ form of the LGR collocation method to discretize the multiple-interval form of the trajectory optimization optimal control problem \cite{Garg_2010,Garg_2011,GargPatterson_2011,Darby_2011,DarbyHager_2011,Patterson_2015,Liu_2018}. The state is approximated using a basis of Lagrange polynomials $\ell_j^{(k)}(\tau)$ \cite{Dennis_2019}
 \begin{equation}
     \m{x}(\tau) \approx \m{X}^{(k)}(\tau) = \sum_{j=1}^{N_k + 1} \m{X}_j^{(k)} \ell_j^{(k)} (\tau),
 \end{equation}
 the derivative of the state with respect to $\tau$ is then
 \begin{equation}
     \frac{d \m{x}^{(k)}(\tau)}{d\tau} \approx \frac{d\m{X}^{(k)}(\tau)}{d\tau} = \sum_{j=1}^{N_k + 1} \m{X}_j^{(k)} \frac{d\ell_j^{(k)} (\tau)}{d\tau},
 \end{equation}
 where
 \begin{equation}
     \ell_j^{(k)}(\tau) = \prod_{l=1,l \neq j}^{N_k + 1} \frac{\tau - \tau_l^{(k)}}{\tau_j^{(k)} - \tau_l^{(k)}}.
 \end{equation}
 The mesh domain is again defined on $\tau \in [-1,+1]$ where $\tau^{(k)}=(\tau_1^{(k)},...,\tau_{N_k}^{(k)})$ are the LGR collocation points on the $k^{th}$ mesh interval and $\tau_{N_k+1}^{(k)} = T_k$ is a noncollocated point. The problem can then be converted to a nonlinear programming problem (NLP) by writing the cost in terms of the $N_k$ collocated points as
 \begin{equation}
     J \approx \C{M}(\m{X}_1^{(1)},t_0,\m{X}_{N_K+1}^{(K)},t_f)
     +\frac{t_f-t_0}{2} \sum_{k=1}^K \sum_{j=1}^{N_k} w_j^{(k)} \C{L}(\m{X}_j^{(k)},\m{U}_j^{(k)},t(\tau_j^{(k)},t_0,t_f)),
 \end{equation}
 where the running cost is approximated using an $N_k$-point LGR quadrature such that $w_j^{(k)}$ are the $N_k$ LGR weights in each mesh interval. The NLP is then subject to the dynamic constraints
 \begin{equation}
     \sum_{j=1}^{N_k+1} D_{ij}^{(k)} \m{X}_j^{(k)} - \frac{t_f - t_0}{2} \m{f}(\m{X}_i^{(k)},\m{U}_i^{(k)},t(\tau_i^{(k)},t_0,t_f)) = \m{0},
 \end{equation}
 the boundary conditions
 \begin{equation}
     \m{b}_{\min} \leq \m{b}(\m{X}_1^{(1)},t_0,\m{X}_{N_K+1}^{(K)},t_f) \leq \m{b}_{\max},
 \end{equation}
 and any path constraints
 \begin{equation}
     \m{c}_{\min} \leq \m{c}(\m{X}_i^{(k)},\m{U}_i^{(k)},t(\tau_i^{(k)},t_0,t_f)) \leq \m{c}_{\max}.
 \end{equation}
 The elements of the LGR differentiation matrix of size $N_k \times (N_k + 1)$ are denoted $D_{ij}^{(k)}$ where $D_{ij}^{(k)} = d\ell_j^{(k)} (\tau_{i}^{(k)})/
 d\tau$ for $(i=1,...,N_k,j=1,...,N_k+1)$.
 
\section{Desensitized Optimal Control}
The objective of desensitized optimal control is to determine the state and control that minimize some performance index along with sensitivities of a user-specified function to state perturbations while satisfying dynamic constraints, boundary conditions, and any path constraints. The user-specified function acts as a penalty term in the cost that quantifies the influence of state perturbations encountered anywhere along the trajectory on the final state. In the context of this research, the state perturbations are a result of parametric uncertainties in the dynamic model. In general, desensitized optimal control involves the introduction of sensitivities, $\m{S}(t)$, in either matrix or function form, as states to the original problem formulation. A sensitivity matrix, as formulated by Ref.~\cite{seywald2003desensitized}, allows for consideration of uncertainties with respect to time varying parameters. As a result, evaluating the cost functional requires propagating $(n \times m)^2 + n + m$ states where $n$ is the number of states in the original optimal control problem formulation and $m$ is the number of parameters with uncertainties. Alternatively, the sensitivity function derived in Ref.~\cite{Makkapati_2018} assumes the parameters of interest to be constant and therefore requires only $n \times m$ states to be propagated. While different approaches have been conceived for introducing sensitivities \cite{seywald2003desensitized,Makkapati_2018}, in this research the approach developed in Ref.~\cite{Makkapati_2018} is employed because it reduces dimensionality compared with the approach developed in Ref.~\cite{seywald2003desensitized}. 

Consider a desensitized optimal control problem with dynamics of the form
\begin{equation}
    \dot{\m{x}}=\m{f}(\m{x},\m{p},\m{u},t),
\end{equation}
which are assumed to be continuous in $(\m{x},\m{p},\m{u},t)$ and continuously differentiable with respect to the state, $\m{x}$, and the nominal parameter values, $\m{p}$. Now suppose the solution to those dynamics is given as 
\begin{equation}
    \m{x}(\m{p},t) = \m{x}_0 + \int_{t_0}^{t_f} \m{f}(\m{x}(\m{p},\tau),\m{p},\m{u}(\tau),\tau)d\tau,
\end{equation}
where the initial condition on the original state vector is known. Let the partial derivative of the state with respect to the parameters now be taken as a function of the elapsed time $t$ and parameters 
\begin{equation}
    \frac{\partial \m{x}}{\partial \m{p}}(\m{p},t) = \int_{t_0}^{t_f} \Biggr[ \frac{\partial \m{f}}{\partial \m{x}}(\m{x}(\m{p},\tau),\m{p},\m{u}(\tau),\tau)\frac{\partial \m{x}}{\partial \m{p}}(\m{p},\tau) 
    + \frac{\partial \m{f}}{\partial \m{p}}(\m{x}(\m{p},\tau),\m{p},\m{u}(\tau),\tau) \Biggr] d\tau.
\end{equation}
Taking the derivative with respect to the elapsed time yields the following expression for the sensitivity dynamics:
\begin{equation}
    \frac{d}{dt} \Biggr[ \frac{\partial \m{x}}{\partial \m{p}} (\m{p},t) \Biggr] = \frac{\partial \m{f}}{\partial \m{x}}(\m{x}(t),\m{p},\m{u}(t),t) \frac{\partial \m{x}}{\partial \m{p}}(\m{p},t)
    + \frac{\partial \m{f}}{\partial \m{p}} (\m{x}(t),\m{p},\m{u}(t),t),
\end{equation}
where the sensitivity is now defined as the change in the state with respect to nominal parameter values
\begin{equation}
    \m{S}(t) = \frac{\partial \m{x}}{\partial \m{p}}(\m{p},t),
\end{equation}
with $n \times m$ elements. The sensitivity dynamics are then
\begin{equation}
    \dot{\m{S}}(t) = \frac{d}{dt} \Biggr[ \frac{\partial \m{x}}{\partial \m{p}} (\m{p},t) \Biggr],
\end{equation}
where the initial condition on the sensitivity matrix is simply the zero matrix since the initial state is fixed. The augmented cost is now a function of the original cost and the sensitivity function at the final time subject to some user-defined weighting term, $Q(t)$, also evaluated at the final time
\begin{equation}
    J_A = J + \int_{t_0}^{t_f} ||\m{S}(t_f)||_{Q(t_f)}^2,
\end{equation}
where, again, the sensitivities are a function of the propagated state under nominal parameter conditions and the nominal parameters are assumed to be constant to improve computational efficiency.

\section{Desensitized Optimal Guidance\label{sect:desensitized}}
The following framework for desensitized optimal guidance combines the desensitized trajectory optimization method from Ref.~\cite{Makkapati_2021} and mesh remapping guidance method from Ref.~\cite{Dennis_2019}. Consider the aforementioned Bolza optimal control problem. Suppose now that that the state, $\m{x}(\m{p},t)$, is subject to parametric uncertainties where the nominal values of the parameters are denoted by $\m{p}$. As a result, it is desirable to reduce sensitivities to these uncertainties so that the resultant error in the final state is minimized. The sensitivity, $\m{S}(t)$, of the state with respect to perturbations at any point along the trajectory is given as the solution to the ordinary differential equation \cite{Makkapati_2018}
\begin{equation}\label{eq: sensitivityDynamics}
    \dot{\m{S}}(t) = \m{A}(t) \m{S}(t) + \m{B}(t),
\end{equation}
where $\m{S}(t_0) = \m{0}$. Here, $\m{A}(t)$ and $\m{B}(t)$ are defined as

%is the partial derivative of the dynamics with respect to the state, and $\m{B}(t)$ is the partial derivative of the dynamics with respect to the nominal parameters.
\begin{equation}
    \m{A}(t) = \frac{\partial \m{f}}{\partial \m{x}}(\m{x}(\m{p},t),\m{u}(t),\m{p},t),
\end{equation}
\begin{equation}
    \m{B}(t) = \frac{\partial \m{f}}{\partial \m{p}}(\m{x}(\m{p},t),\m{u}(t),\m{p},t).
\end{equation}
To minimize perturbations in the final state, the objective functional must be augmented, $J_A$, to include variations in a user-defined penalty term, $\m{h}=\m{g}(x)$.
\begin{equation}\label{eq: augmentedCost}
    J_A = J + \mathbb{E}\left(||\delta \m{h}(t_f)||_{\m{Q}_f}^2 + \int_{t_0}^{t_f} ||\delta \m{h}(t)||_{\m{Q}(t)}^2 dt \right)
\end{equation}
The penalty error term is then defined as \cite{Makkapati_2021}
\begin{equation}
    \mathbb{E}(||\delta \m{h}||_{\m{Q}}^2) = \text{tr}\m{Q} \mathbb{E}(\delta \m{h} \delta \m{h}^\top) \approx \text{tr} \m{Q} \m{G} \m{S} \m{P} \m{S}^\top \m{G}^\top, 
\end{equation}
where $\m{Q}$ is a user-defined weighting matrix that is positive semi-definite, $\m{G}$ is the jacobian of the penalty term, $\m{h}$, and $\m{S}$ is the sensitivity function with $n \times m$ elements. $\m{P}$ is a user-defined, positive semi-definite covariance matrix which is a function of the nominal values of the parameters. This formulation allows for the user to define the uncertainty in the parameter through $\m{P}$ while adjusting how much it is desired to desensitize the control through $\m{Q}$. If $\m{Q}$ is set to zero at all points in time, then the objective returns to that of the original optimal control problem.

\begin{figure}[h]
\centerline{\includegraphics[scale=0.72]{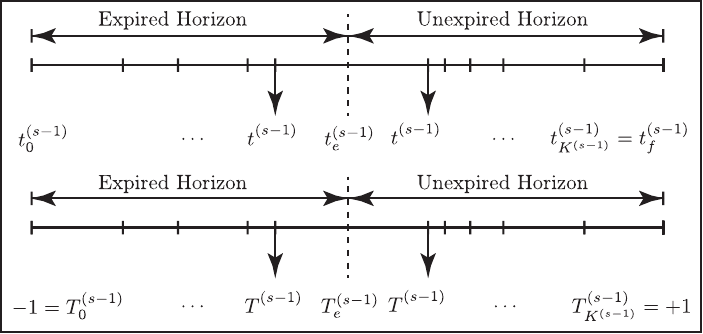}}
\caption{Mesh remapping strategy for guidance cycle $(s-1)$ featuring the expired and unexpired horizon \cite{Dennis_2019}.}
\label{fig:meshRemapping}
\end{figure}

The aforementioned problem formulation determines a desensitized reference trajectory and therefore a control which is less sensitive to perturbations in the state due to parametric uncertainties. Suppose now that it is desired to perform guidance updates such that the desensitized optimal trajectory is recalculated on the remaining horizon to allow for corrections in the desensitized optimal control. Let $s$ denote the current guidance cycle where $s \in [1,2,3,...,S]$ and $D$ denote the duration of the guidance cycle. The current guidance cycle iteration then occurs on the time interval $t \in [t_0^{(s)}, t_e^{(s)}]$ where
\begin{equation}
    t_0^{(s)} = t_0 + s D,
\end{equation}
\begin{equation}
    t_e^{(s)} = t_0 + (s + 1) D.
\end{equation}
At the end of each guidance cycle, the expired horizon is removed and the initial conditions updated before the problem can be resolved for the remainder of the mesh. The terminal time of the previous cycle occurs between two mesh points as seen in Fig. \ref{fig:meshRemapping}. To delete the expired horizon, the desensitized optimal control, state, and sensitivity are interpolated to $t_e^{(s-1)}$ and the initial conditions become
\begin{equation}
    \m{x}_0 = \m{x}(t_0^{(s)}) = \Tilde{\m{x}}(t_e^{(s-1)}),
\end{equation}
\begin{equation}
    \m{S}_0 = \m{S}(t_0^{(s)}) = \m{S}(t_e^{(s-1)}),
\end{equation}
while the final boundary conditions remain the same. The mesh is then remapped such that the first mesh point corresponds to $t_e^{(s-1)}$ and $t \in [t_0^{(s)},t_f]$ maps to $\tau \in [-1, +1]$. Note that $t_e^{(s)}$ is not the same as $t_f$. The former is the final time for the simulated dynamics of the relevant guidance cycle while the latter is the terminal boundary condition on the time for the desensitized optimal control problem.  

To reiterate, the desensitized optimal control problem is first solved using the nominal parameter values to obtain a desensitized optimal control for the entirety of the trajectory. The dynamics are then simulated for the first guidance cycle using the perturbed parameter values and the reference desensitized optimal control. At the end of the guidance cycle, the final value of the simulated dynamics is used as the initial state at the start of the next guidance cycle. The expired horizon is removed and the mesh is remapped after interpolating the desensitized reference state, control, and sensitivity to the current time. The problem is then resolved on the remaining horizon to obtain a new desensitized optimal control used to simulate the dynamics for the next cycle. This process repeats until the end of the horizon with the goal of minimizing sensitivities to parametric uncertainties while allowing for corrections in the control to reduce perturbations in the final state. A simple numerical example is introduced in the next section to demonstrate this guidance method.

\begin{figure}[h]
\centerline{\includegraphics[scale=0.72]{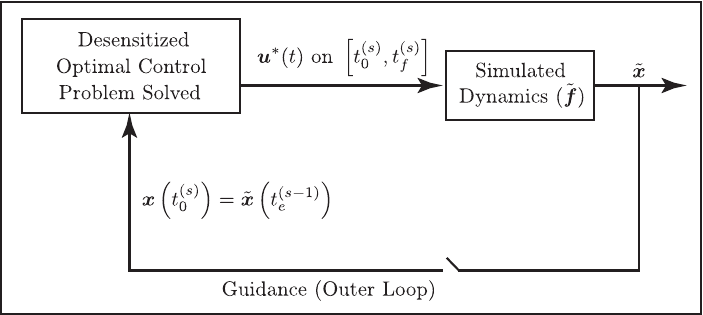}}
\caption{Block diagram of desensitized optimal guidance method where $s$ is the current guidance cycle iteration.}
\label{fig:DOG}
\end{figure}

\section{Example}

As an example of applying the approach described in Section~\ref{sect:desensitized}, consider the following optimal control problem: % Determine the state, $x(t)$, and control, $u(t)$, that 
\begin{equation}\label{example}
   \min~ J = \frac{1}{2}\int_{t_0}^{t_f} \left(x^2(t) + u^2(t) \right)dt ~~\textrm{subject to}~~\left\{\begin{array}{lcl} \dot{x} & = & - \alpha^2 x^3 + \alpha u, \\ (x(t_0),x(t_f)) & = & (1.5,1.0), \\ (t_0,t_f) & = & (0,50). \end{array} \right.
\end{equation}
% subject to the dynamic constraint
% \begin{equation}\label{example-dynamics}
%   \dot{x} = - \alpha^2 x^3 + \alpha u,
% \end{equation}
% and boundary conditions
% \begin{equation}\label{example-bcs-and-parameters}
%     \begin{array}{lcl}
%         x(t_0) & = & x_0 = 1.5,\\
%         x(t_f) & = & x_f = 1.0,\\
%         t_0 & = & 0,\\
%         t_f & = & 50 \\
%         \alpha & = & 2.
%     \end{array}
% \end{equation}
Suppose now that it is desired to design both a reference solution and a guidance solution that reduces sensitivities in the terminal state with respect to perturbations in the parameter $\alpha$.  A desensitized optimal control problem that meets this aforementioned objective is then given as follows.  The objective functional in Eq.~\eqref{example} is then modified to be
\begin{equation}\label{desensitized-objective}
    J_A = \frac{1}{2}\int_{t_0}^{t_f} \left(x^2(t) + u^2(t) \right)dt + \mathbb{E}\left(||\delta h(t_f)||_{Q_f}^2 \right),
\end{equation}
% subject to the dynamics and boundary conditions given in Eq.~\eqref{example} with $S(t_0)=0$.
Because the objective in desensitizing the control is to minimize the final state error, the penalty term $h$ is chosen as
\begin{equation}
    h = (x_f),
\end{equation}
and $G=1$ since the problem contains only one state. The parameter covariance is a function of the nominal value and standard deviation of the parameter of interest. Let $q$ denote a user-defined percentage of the nominal value of the parameter. The standard deviation is then set such that $\sigma = q \alpha$. The covariance is then $P = \sigma^2$. The weighting function, $Q(t)$, is a constant weighting parameter defined as $Q(t) = \beta$ for each simulation where $\beta$ can be altered depending on the magnitude of desensitization desired. When $\beta$ is assigned a value of zero, the problem returns to a standard optimal control problem. The larger the value of $\beta$, the more desensitization there is. 

\begin{figure}[h]
  \subfigure[$(\sigma,\beta)=(0.01\alpha,5)$.\label{fig:errorPlots-B=5 f=0.01}]{\includegraphics[scale=0.4]{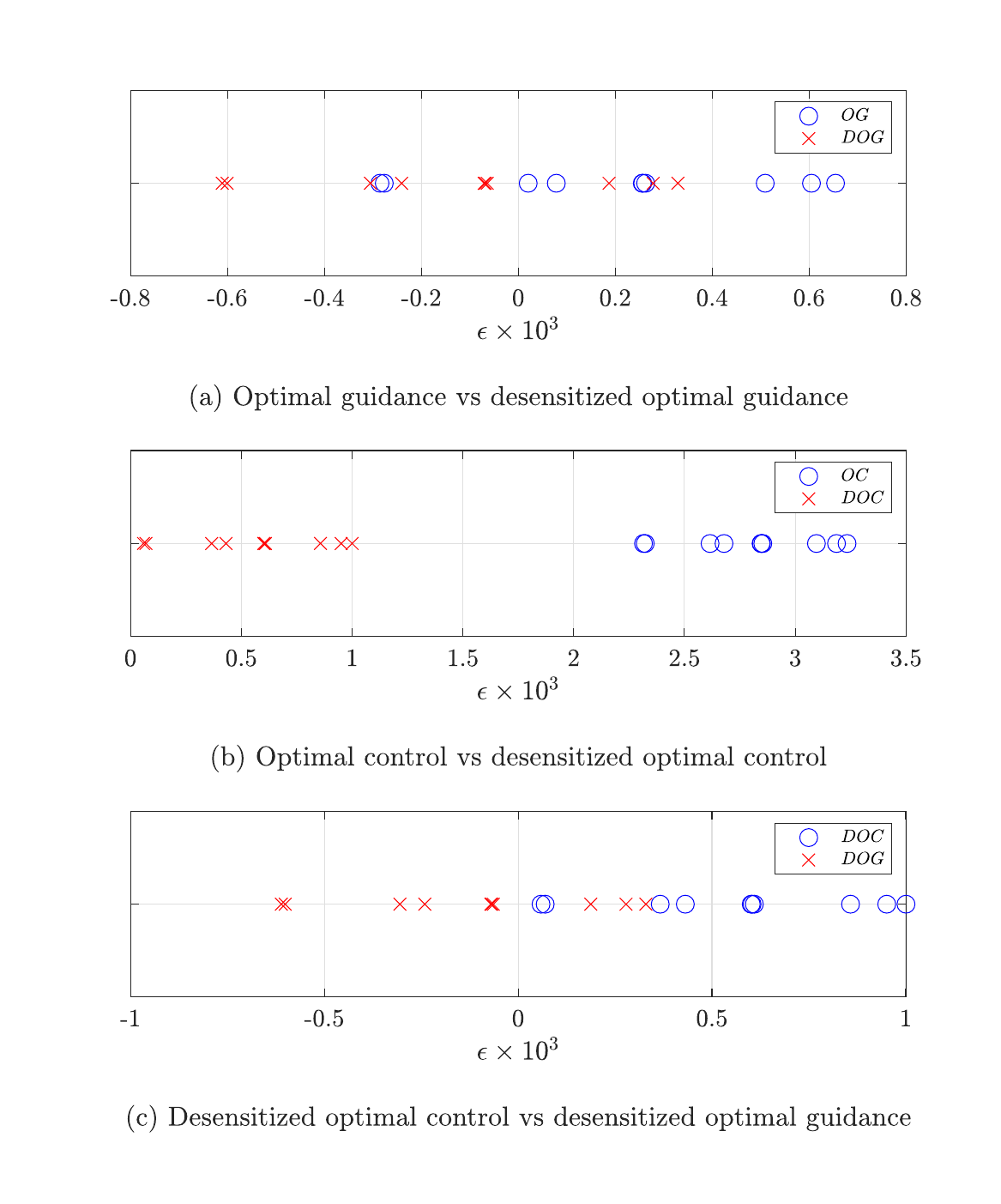}}
  \subfigure[$(\sigma,\beta)= (0.01\alpha,10)$.\label{fig:errorPlots-B=10 f=0.01}]{\includegraphics[scale=0.4]{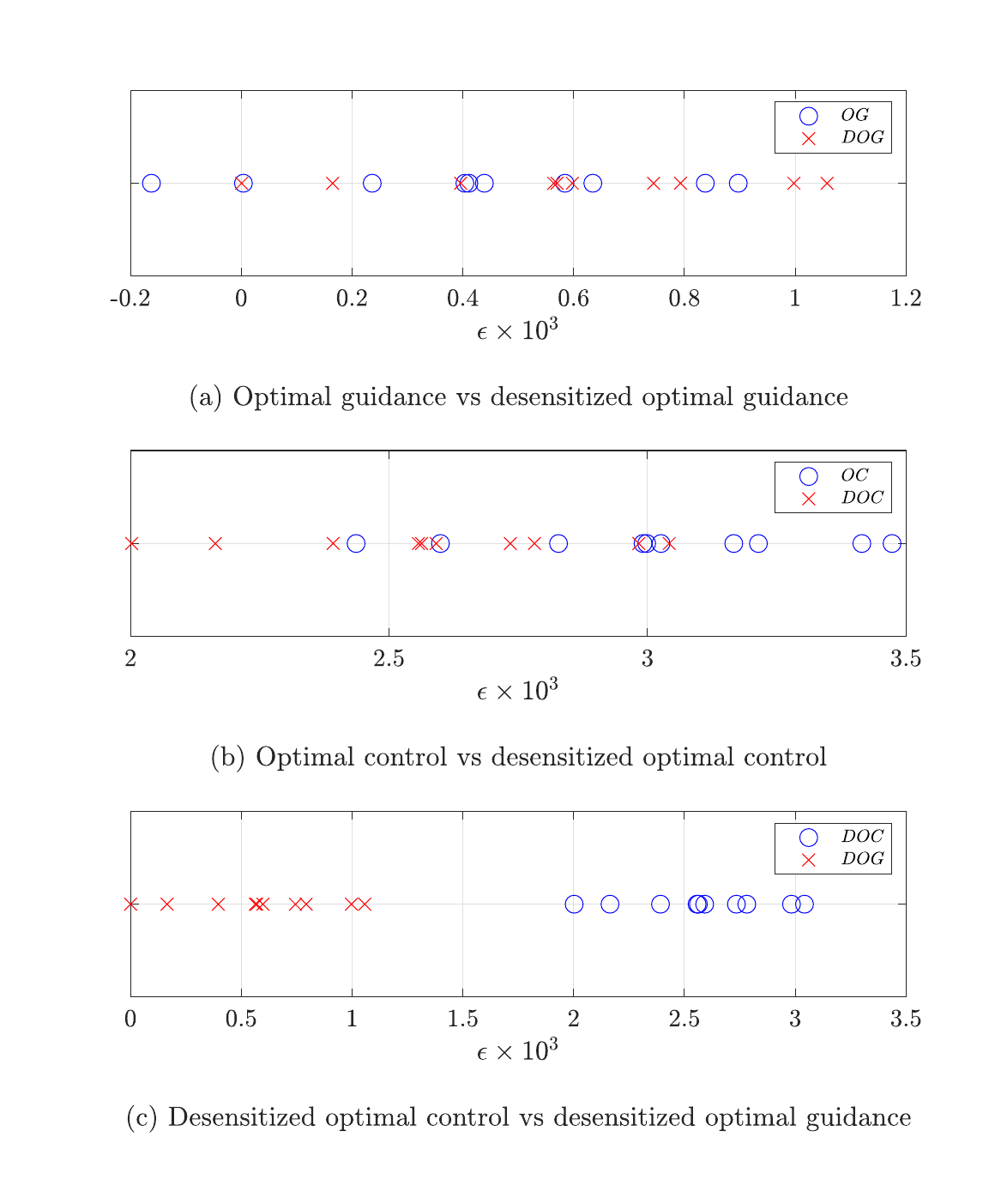}}

  \subfigure[$(\sigma,\beta)= (0.02\alpha,5)$. \label{fig:errorPlots-B=5 f=0.02}]{\includegraphics[scale=0.4]{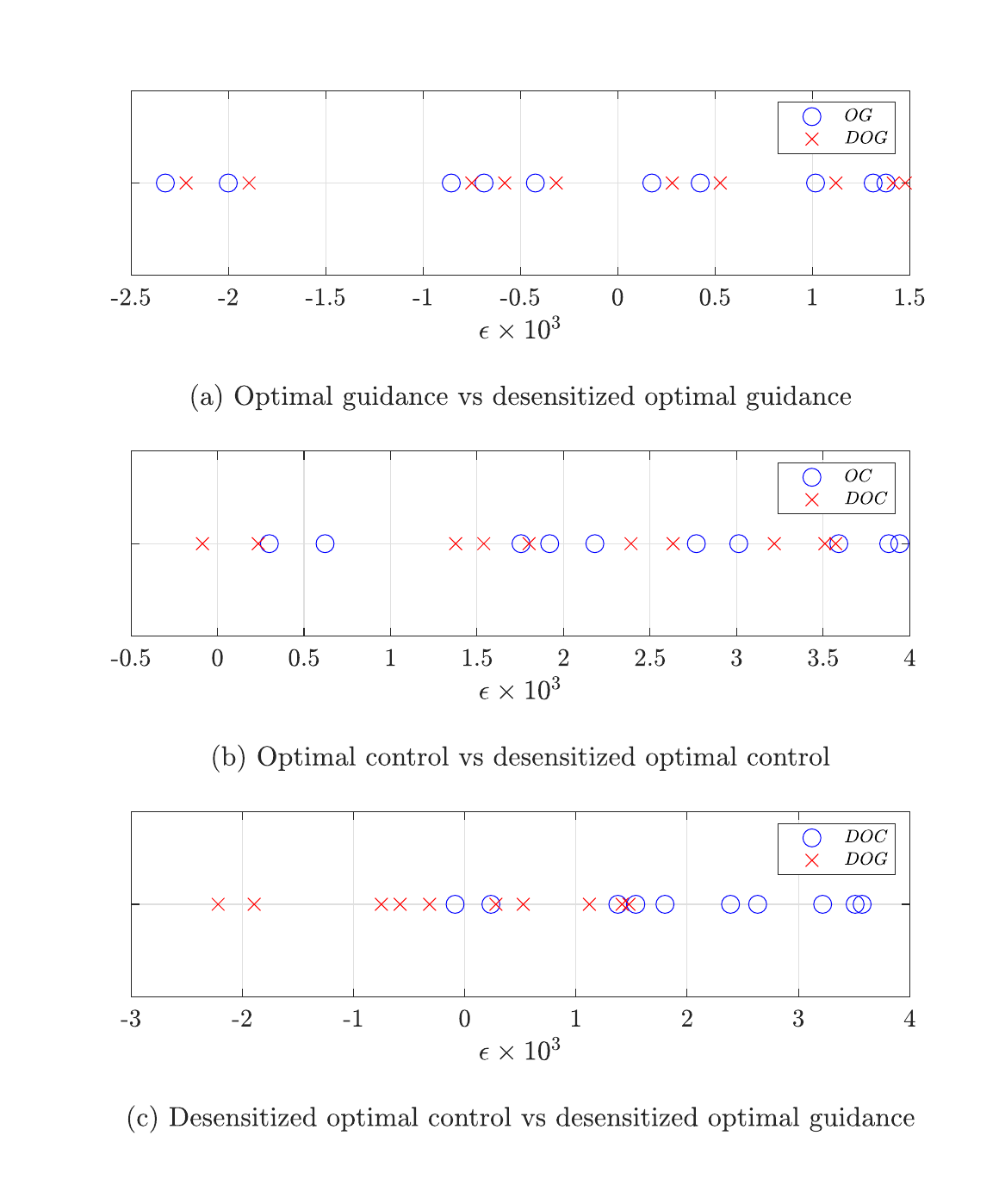}}
  \subfigure[$(\sigma,\beta)= (0.02\alpha,10)$. \label{fig:errorPlots-B=10 f=0.02}]{\includegraphics[scale=0.4]{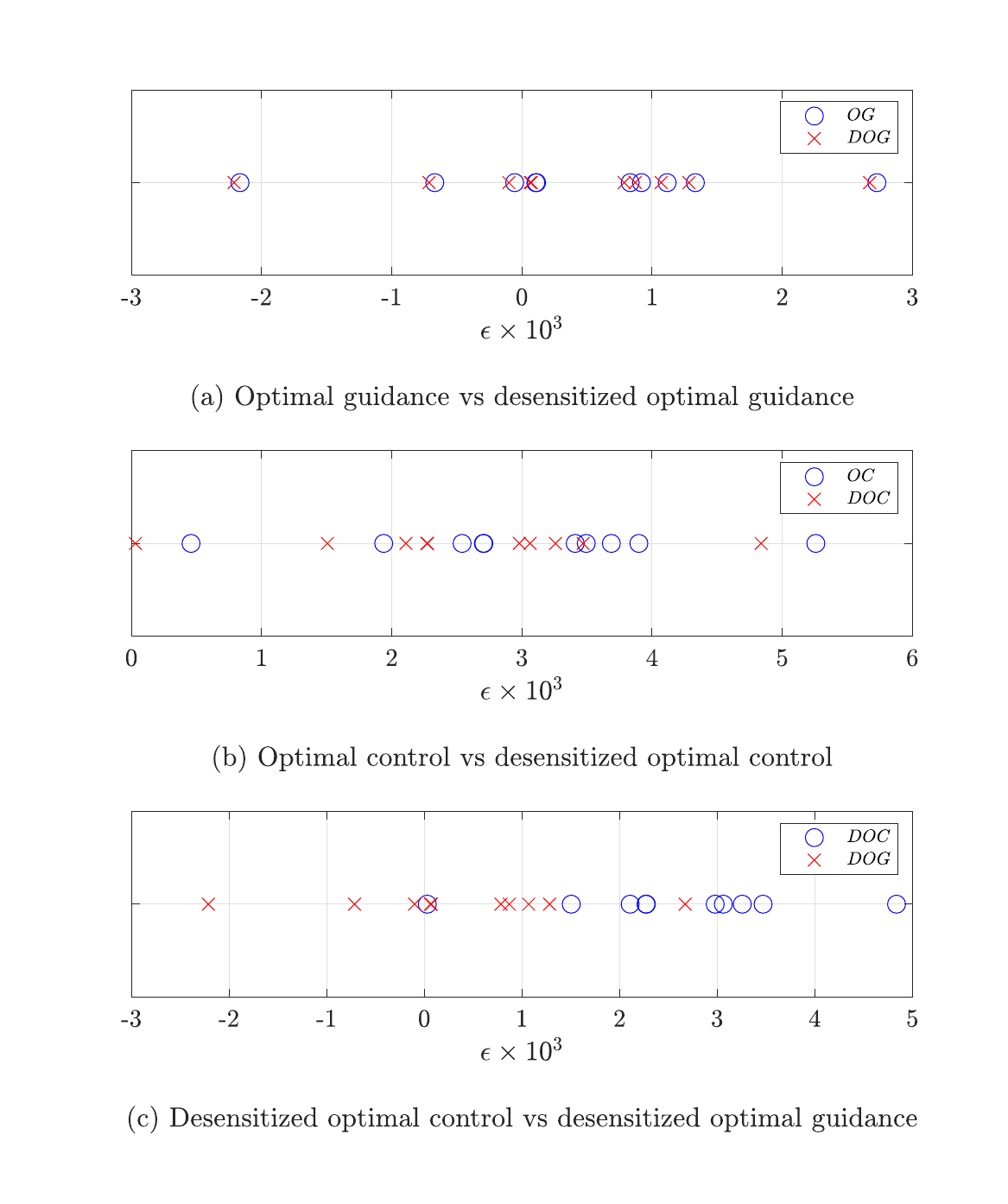}}

  \caption{Final state error results for $(\sigma,\beta) = \left\{(0.01\alpha,5),(0.01\alpha,10),(0.02\alpha,5),(0.02\alpha,10)\right\}$.\label{error-plots}}
\end{figure}

To evaluate the effectiveness of the desensitized optimal guidance method, the following four sets of 100 Monte Carlo simulations were performed for each combination of the weighting parameter and standard deviation: $(\alpha,\beta)=\left\{(0.01\alpha,5),(0.01\alpha,10),(0.02\alpha,5),(0.02\alpha,10)\right\}$.  It is noted, however, that only ten simulations per case are presented for clarity. Furthermore, for this study, the parameter $\alpha$ was sampled from the Gaussian distribution
\begin{equation}
    \Tilde{\alpha} \sim N(\alpha,P).
\end{equation}
For the full set of Monte Carlo simulations, the deviation in the final state relative to the reference, denoted $\epsilon$ (which implies that $\epsilon$ can be negative), is then shown in Fig.~\ref{error-plots}.  Each point shown for each case in Fig.~\ref{error-plots} represents the deviation in the final state from that of the reference trajectory for a single simulation. For each case, the desensitized optimal guidance scheme was compared to the methods presented in Ref.~\cite{Makkapati_2021} and \cite{Dennis_2019}. Figures \ref{fig:errorPlots-B=5 f=0.01}--\ref{fig:errorPlots-B=10 f=0.02} present the results for the aforementioned four cases of $(\alpha,\beta)$, respectively.  Results are compared for standard optimal control with guidance (OG), desensitized optimal control with guidance (DOG), optimal control without guidance (OC), and desensitized optimal control (DOC). Note that Fig.~\ref{fig:errorPlots-B=10 f=0.01} serves primarily to show that desensitized optimal control works for trajectory optimization. For methods with guidance updates, the guidance cycle duration was set to four seconds with a total of 12 cycles performed. 

\begin{figure}[h]
\centerline{\includegraphics[scale=0.47]{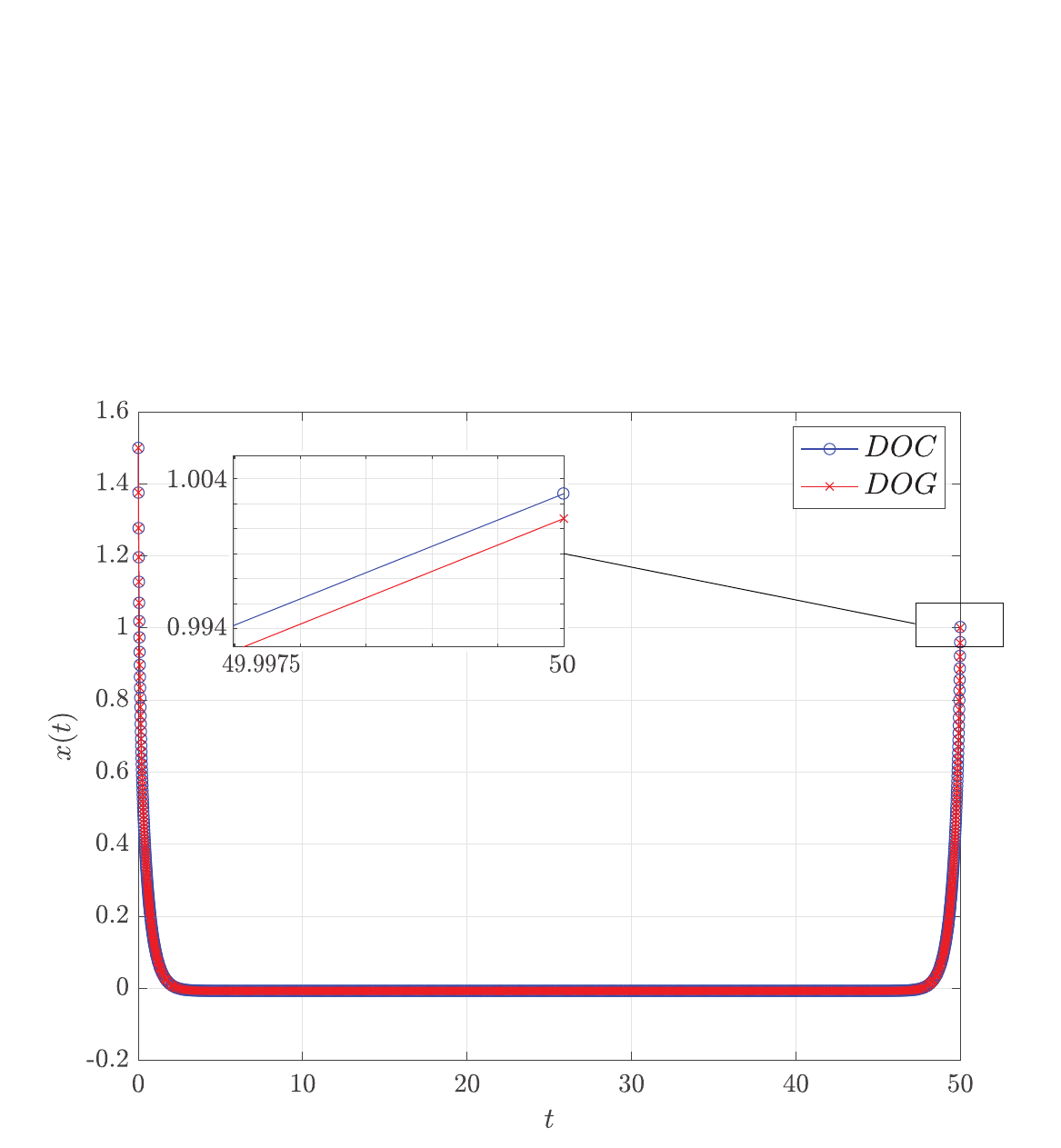}}
\caption{DOC state vs DOG state for $\beta=10$ and $\sigma=0.01\alpha$ for $\tilde{\alpha}=2.0178$.}
\label{fig:DOGvsDOC state for B=10 f=0.01}
\end{figure}

The numeric results displayed demonstrate that desensitized optimal control with guidance updates performs overall better than both desensitized optimal control without guidance updates and standard optimal control with guidance updates. This can most clearly be seen in Fig.~\ref{fig:errorPlots-B=5 f=0.01} where the center of the Monte Carlo distribution for DOG is closer to zero than that of OG and DOC. Thus, it is apparent that combining the methods of Ref.~\cite{Makkapati_2021} and \cite{Dennis_2019} reduces terminal state error as compared with the results of each individual method. It should be noted that the results are highly dependent on the user-defined weighting parameters and relative magnitude of the parametric uncertainties. If $Q(t)$ is weighted too highly and the perturbations in the parameters are not large enough, then OG will perform better than DOG as seen in Fig.~\ref{fig:errorPlots B=10 f=0.01}(a). Thus, the correct balance of desensitization is required to correctly address perturbations while still producing the optimal solution. Finally, the trajectories obtained for a single Monte Carlo simulation are shown in Fig.~\ref{fig:DOGvsDOC state for B=10 f=0.01} where an enlarged view of the final state is included with and without guidance updates. 

\section{Conclusions}
A guidance method has been developed using desensitized optimal control in conjunction with guidance updates. The optimal control was desensitized to state perturbations resulting from parametric uncertainties to promote robustness. For each guidance update performed, the expired horizon was removed and the mesh was remapped on the remaining horizon. The desensitized optimal control was then resolved to allow for control corrections. This approach was applied to a simple example, and numerical results demonstrated that combining desensitized optimal control with guidance updates produces results where the center of distribution for final state error is closer to zero than either that of desensitized optimal control without guidance updates or using a previously developed optimal guidance method. Future work will apply this method to a more complicated problem to further explore its efficiency in reducing state perturbations for high performance vehicles.

\clearpage
\bibliographystyle{ieeetr}
% \input{WinklerRaoCDC-ArXiV.bbl}

% \bibliography{references}

\end{document}